\documentclass{article}

\usepackage{amssymb}
\usepackage{amsmath}

\usepackage[english,francais]{babel}

\newtheorem{theorem}{Theorem}[section]

\newtheorem{e-proposition}[theorem]{Proposition}

\newtheorem{e-definition}[theorem]{Definition\rm}

\newcommand{\Z}{{\mathbb{Z}}}   
\newcommand{\Q}{{\mathbb{Q}}}   
\newcommand{\M}{{\mathcal{M}^{0}}}   
\newcommand{\U}{{\mathcal{U}}}   
\newcommand{\vp}{\vspace{0,5cm}}   

\newcommand{\A}{{\mathcal{A}}}
\newcommand{\I}{{\mathfrak{I}}}

\newcommand{\s}{{\mathcal{S}}}  
\newcommand{\p}{{\mathcal{P}}}  
\def\og{\leavevmode\raise.3ex\hbox{$\scriptscriptstyle\langle\!\langle$~}}
\def\fg{\leavevmode\raise.3ex\hbox{~$\!\scriptscriptstyle\,\rangle\!\rangle$}}

\begin{document}
\centerline{}


\selectlanguage{english}
\begin{center}
\Large \textbf{On the Classification of Finite Semigroups and \emph{RA}-loops
  with the Hyperbolic Property }
\normalsize $$\textrm{S. O. Juriaans}^{a} \qquad \textrm{A. C. Souza
  Filho}^{b}$$
\end{center}


\selectlanguage{english}

\begin{center}
{\it Instituto de Matem\'atica e Estat\'\i stica,
 Universidade de S\~ao Paulo,
Caixa Postal 66281, S\~ao Paulo, CEP 05315-970 - Brazil}
\end{center}
\begin{center}
\small
email adresses:
$^{a}$ostanley@ime.usp.br $\quad$ $^{b}$calixto@ime.usp.br
\end{center}
\normalsize



\begin{abstract}

 We classify the finite semigroups $S$, for
which all $\Z$-orders $\Gamma$ of $\Q S$, the unit group $\U(\Gamma)$
is hyperbolic. We also  classify the $RA$-loops $L$ for which the
unit loop of its integral loop ring does not contain any free abelian subgroup
of rank two.
\end{abstract}







\section{Introduction}
Initially we classify finite dimensional algebras $\A$
over $\Q$ such that if $\Gamma \subset \A$ is a $\Z$-order then  $\U(\Gamma)$
is hyperbolic. If $\A$ is such an algebra, we say that $\A$ has the hyperbolic property.

\vp

In \cite{gdf}, are classified the semigroups $\Sigma$ with $\U (\Z \Sigma)$
finite. Therefore, for this class of semigroups the algebra $\Q
\Sigma$ has the hyperbolic property. 

\vp

First we classify the semisimple algebras $\Q S$ which are nilpotent free. If
$\Q S$ has nilpotent elements there are two possibilities: either $S$ contains
nilpotent elements, or not. In the latter case $S$ is a disjoint
union of groups of certain types. In the former,
$S$ is a union of groups and a subsemigroup of order five.
We also give the structure of $S$ when $\Q S$ is  non-semisimple and has the
hyperbolic property.
For the proofs of the results of section two see \cite{hsa}.

 \vp

 In the last section we classify the $RA$-loops $L$, such that, $\Z^2 {\not \hookrightarrow} \U (\Z L)$.

\vp



\vp
\section{Semigroup Algebras}

We will consider $\A$ a unitary finitely generated $\Q$-algebra and 
denote by  $\s (\A)$, respectively  $J(\A)$, the semisimple
subalgebra, respectively the Jacobson radical, of $\A$ and by 
$E(\A)=\{\ E_{1}, \cdots, E_{N}\}, N \in \Z^{+}$, the set of the
central primitive idempotents of the semisimple algebra $\s(\A)$. A
classical result of Wedderburn-Mal\^cev states that 
$$\A \cong \s (\A) \oplus J(\A),$$   as a vector space. As a result,
we have that  $\A$ is a artinian algebra and thus its radical  is a 
nilpotent ideal. We denote
$T_2 (\Q ):=\left ( \begin{array}{ll}\Q & \Q \\ 0 & \Q \end{array}
\right )$ the $2 \times 2$ upper triangular matrices over $\Q$, with
the usual matrix multiplication.

\begin{e-definition}
Let $\A$ be a finite dimensional algebra over $\Q$ and $\Gamma$ be a $\Z$-order of
$\A$. If  $$\Z^{2} {\not \hookrightarrow} \U (\Gamma),$$ we say $\A$ has the
hyperbolic property.
\end{e-definition}

\vp

\begin{theorem}
Let $\A$ be a finite dimensional $\Q$-algebra. If $\A_i$ is a simple
epimorphic image of $\A$, denote by $F_i$ a maximal subfield of
$\A_i$ and $\Gamma_i \subset \A_i$ a $\Z$-order.
 The following conditions hold:
\begin{enumerate}
\item The algebra $\A$ has the hyperbolic property, it is semisimple and it
  has no nilpotent element  if, and only if,  $$\A =\oplus \A_{i},$$ whereof $\A_{i}$ is a
  division ring and there exists at most one index $i_{0}$ such that $\U(\Gamma_{i_{0}})$ is hyperbolic and  infinite.
\item The algebra $\A$ has the hyperbolic property and it is semisimple with
  nilpotent elements if, and only if,
$$\A =(\oplus \A_{i})\oplus M_2(\Q).$$
\item The algebra $\A$ has the hyperbolic property and it is non semisimple with
  central radical if, and only if, $$\A =(\oplus \A_{i})\oplus J.$$
\item The algebra $\A$ has the hyperbolic property and it is non semisimple with
  non central radical if, and only if, $$\A =(\oplus \A_{i})\oplus
 T_{2}(\Q).$$
\end{enumerate}
 For each item above, $F_i$ is an imaginary quadratic field and
  $\A_{i}$ is either  an imaginary quadratic field or a totally
  definite quaternion algebra. Furthermore, every  simple
epimorphic image of $\A$ in the direct sum is an
  ideal of $\A$.
\end{theorem}

\vp

In what follows, $S$ denotes a finite semigroup, $\Q S$ denotes a
unitary semigroup algebra over $\Q$, $\M(G;n,n;P)$ denotes the Rees semigroup with structural group
$G$ and $P$ denotes an $n \times n$ sandwich matrix.

\vp

\begin{theorem} \label{tsmp}
The algebra $\Q S$ has no nilpotent element and it has the hyperbolic property if, and
only if, $S$ is an inverse semi-group and it admits a principal series, whose principal factors are
isomorphic to groups $G$ and at most a unique $K$, listed below:
\begin{enumerate}
\item $G$ is an abelian group of exponent dividing $4$ or $6$;
\item $G$ is a hamiltonian $2$-group;
\item $K \in \{C_{5},C_{8},C_{12}\}$.
\end{enumerate}
\end{theorem}
 
\vp

\begin{theorem}
Let $\Q S$ be an algebra with nilpotent elements. The algebra $\Q S$ is
semisimple and it has the hyperbolic property if, and only if, $S$ admits a
principal series whose principal factors are isomorphic to groups $G$ and a
unique semigroup $K$, listed below:
\begin{enumerate}
\item  $G$ is an abelian group of exponent dividing $4$ or $6$.
\item $G$ is a hamiltonian $2$-group.
\item $K$ is a group of the set $\{S_{3},D_{4},Q_{12},C_{4}\rtimes C_{4}\}$.

\item $K$ is one of the Rees semigroups: $$\M(\{1\};2,2;I_{d})=M \quad \quad \textrm{ or }\quad \quad\M(\{1\};2,2; \left(\begin{array}{ll}1&1\\0&1\end{array}\right))=M_{12},$$ which is an ideal of  $S$.
\end{enumerate}
In particular, $S$ is the disjoint union of the groups $G$ and the semigroup $K$.

\end{theorem}

\vp

\begin{theorem}
The algebra $\Q S$ is non semisimple and it has the hyperbolic property
if, and only if, there exists a unique nilpotent element $j_{0} \in
S$, such that, the subsemigroup $\I= \{\theta, j_{0}\}$ is an ideal
of $S$, and $S \setminus \I$ admits a principal series whose
principal factors are isomorphic to abelian groups of exponent
dividing $4$ or $6$, or a hamiltonian $2$-group. In particular $S/
\I$ is the disjoint union of its maximal subgroups such that if
$e_{1} \in G_{1}$, and $e_{N} \in G_{N}$ are the respective group identity element, then $e_1 j_0 =j_0 e_N=j_0$. Writing
$$\begin{array}{l}e_1=\sum E_{1_{l}}+  E_{1}+\lambda j_{0}, \lambda\in \Q\\
                  e_N=\sum E_{N_{l}}+  E_{N}+\mu j_{0}, \mu \in \Q \end{array}$$
then only one of the following holds:
\begin{enumerate}
\item $$e_1 e_N =0 \Leftrightarrow e_{N}e_{1}=0 \textrm{ and } \lambda+\mu=0;$$
 $T_2\cong\{ e_{1},e_{N},j_0,\theta\}$ is such that $\Q T_2 \cong T_2(\Q)$.

\item $$\textrm{If } e_N e_1 \neq 0 \textrm{ then } e_1 e_N = e_N e_1=:e_3
  \textrm{ and } \lambda+\mu=0;$$
 $T'_2=\{e_1,e_2,e_3,j_0,\theta\}$ is a subsemigroup of $S$ and $\Q T_2' \cong \Q \oplus \Q \oplus T_{2}(\Q).$
\item $$e_N e_1 =0 \Leftrightarrow e_{1}e_{N}=j_0   \Leftrightarrow \lambda+\mu=1; $$
 $\hat{T}_{2}=\{ e_{1},e_{N},j_0,\theta\}$, and $\Q
\hat{T}_2 \cong T_2(\Q)$.
\end{enumerate}
The semigroups $T_2,T'_2$ and $\hat{T}_2$ are non isomorphic.
\end{theorem}

 \vp

\section{The hiperbolicity of the $RA$-loop loop units}

\vp

In this section we classify the $RA$-loops $L$ such that
$\Z^{2} {\not \hookrightarrow} \U(\Z L)$, the loop of units of $\Z
L$. A loop $L$ is a nonempty set, with a closed binary
operation $\cdot$, such that the equation $a \cdot b=c$ has a
unique $b \in L$ when $a,c \in L$ are known, and a unique solution $a \in L$
when $b,c \in L$ are know, and with a two-side identity
element $1$. Denoting by
$[x,y,z]\dot{=}(xy)z-x(yz)$, recall that a ring $A$ is alternative if
$[x,x,y]=[y,x,x]=0$, for every $x,y \in A$. An $RA$-loop is a loop whose loop ring $R L$ over some
commutative, associative and unitary ring $R$ of characteristic
not equal to $2$ is alternative, but not associative. The basic reference is
\cite{gjp}.

For a theoretical group property $\p$, a group $G$ is virtually $\p$
if  it has  a subgroup of  finite index, $H$ say,  with property $\p$.

\begin{theorem}[\cite{thss}, Theorem 3.3.6]\label{kpm}
Let $L$ be a $RA$-loop.   $\U(\Z L)$ has the hyperbolic property if,
 and only if, $L$ is a finite loop or a loop whose center is virtually
 cyclic, the torsion subloop $T(L)$ of $L$ is such that, if $T(L)$ is a group, then it is an abelian group of exponent dividing $4$
 or $6$ or a hamiltonian $2$-group whose subgroups are all normal in $L$  and if $\ T(L)$ is a loop then it is a
 hamiltonian Moufang $2$-loop whose subloops are all normal in $L$. In this conditions we also have that 
 $\U_{1}(\Z L)=L$.
\end{theorem}

\vp

\section*{Acknowledgements}
This work is part of the second authors Ph.D
thesis, \cite{thss}. He would like to thank his thesis supervisor Prof. Dr. Stanley
Orlando Juriaans for his guidance during this work.

\end{document}